                    	\def\version{29 October, 2013}	                          %
\documentclass[11pt, reqno]{amsart}
\usepackage[centertags]{amsmath}
\usepackage{amsthm, a4, latexsym, amssymb}
\usepackage[latin1]{inputenc}
\usepackage{color}
\usepackage{graphicx}
\usepackage{dsfont}

\def\@rmrk#1#2{\refstepcounter
    {#1}\@ifnextchar[{\@yrmrk{#1}{#2}}{\@xrmrk{#1}{#2}}}

\makeatletter\@addtoreset{equation}{section}\makeatother

\newfont{\bfit}{cmbxti10 scaled 2000}
\newfont{\biggi}{cmr12 scaled 2000}
\newtheorem{step}{STEP}

\newcommand{\bes}{\begin{step}}
\newcommand{\es}{\end{step}}
\newcommand{\X}{\mathcal{X}}

 \newcommand{\eps}{\varepsilon}

 \newcommand{\supp}{{\rm supp}\,}
 \newcommand{\e}{{\rm e}}

 \newcommand{\R}{\mathbb{R}}
 \newcommand{\Z}{\mathbb{Z}}
 \newcommand{\N}{\mathbb{N}}

 \renewcommand{\d}{\,\mathrm{d}}
 
 \newcommand{\prob}{\mathbb{P}}
 \newcommand{\Prob}{{\rm Prob }}
 
 \newcommand{\T}{{\mathcal T}}

 \renewcommand{\P}{\mathbb{P}}
 \newcommand{\E}{\mathbb{E}}
\newcommand{\A}{\mathcal{A}}

 \newcommand{\Dcal}{{\mathcal D}}
 \newcommand{\Acal}{{\mathcal A}}
 \newcommand{\Mcal}{{\mathcal M}}

 \newcommand{\Ical}{{\mathcal I}}

 \newcommand{\Tcal}{{\mathcal T}}
 \newcommand{\Xcal}{{\mathcal X}}
 \newcommand{\Scal}{{\mathcal S}}

 \newcommand{\ssup}[1] {{\scriptscriptstyle{({#1}})}}
\def\1{{\mathchoice {1\mskip-4mu\mathrm l}      % Blackboard bold 1
{1\mskip-4mu\mathrm l}
{1\mskip-4.5mu\mathrm l} {1\mskip-5mu\mathrm l}}}

\newenvironment{proofsect}[1] {\vskip0.1cm\noindent{\rmfamily\bf#1.}}{\qed\vspace{0.15cm}}%{\newline\vspace{0.15cm}}

\newcommand{\ms}{{\Mcal_1^{\ssup {\rm s}}(\T^2)}}
\newcommand{\mms}{{\Mcal_1^{\ssup {\rm s}}((\T\times \mathcal{X})^2)}}

 \newcommand{\kommentar}[1]{}
\renewcommand{\subsection}{\secdef \subsct\sbsect}
\newcommand{\subsct}[2][default]{\refstepcounter{subsection}
\vspace{0.15cm}
{\flushleft\bf \arabic{section}.\arabic{subsection}~\bf #1  }
\nopagebreak\nopagebreak}
\newcommand{\sbsect}[1]{\vspace{0.1cm}\noindent
{\bf #1}\vspace{0.1cm}}

\newtheorem{theorem}{Theorem}[section]

\newtheorem{lemma}[theorem]{Lemma}

\newtheorem{prop}[theorem]{Proposition}
\newtheorem{Assumption}[theorem]{Assumption}

\theoremstyle{definition}

\def\thebibliography#1{\section*{Bibliography}
  \list%
  {\arabic{enumi}.}%                          *** style of reference number ***
    {\settowidth\labelwidth{[#1]}\leftmargin\labelwidth
    \advance\leftmargin\labelsep
    \parsep0pt\itemsep0pt
    \usecounter{enumi}}
    \def\newblock{\hskip .11em plus .33em minus .07em}
    \sloppy                   % \clubpenalty4000\widowpenalty4000
    \sfcode`\.=1000\relax}

\def\P{\prob}

\begin{document}
\title[Moment asymptotics for multitype branching random walks in random environment]{\Large Moment asymptotics\\\medskip for multitype branching random walks\\\medskip in random environment}
\author[Onur G\"un, Wolfgang K\"onig and Ozren Sekulovi\'c]{}
\maketitle
\thispagestyle{empty}
%\vspace{-0.5cm}

\centerline{\sc By  Onur G\"un\footnote{Weierstrass Institute Berlin, Mohrenstr. 39, 10117 Berlin, Germany, {\tt guen@wias-berlin.de} and {\tt koenig@wias-berlin.de}}, Wolfgang K\"onig\footnotemark[1]$^,$\footnote{Institute for Mathematics, TU Berlin, Str.~des 17.~Juni 136, 10623 Berlin, Germany, {\tt koenig@math.tu-berlin.de}}, and Ozren Sekulovi\'c\footnote{University of Montenegro, Cetinjska 2, 81 000 Podgorica, Montenegro, {\tt ozrens@t-com.me}}}
\renewcommand{\thefootnote}{}


\setcounter{footnote}{3}
\renewcommand{\thefootnote}{\arabic {footnote}}
\vspace{-0.5cm}
\centerline{\textit{Weierstrass Institute Berlin, TU Berlin, and University of Montenegro}}
\vspace{0.2cm}

\begin{center}
\version
\end{center}

\begin{quote}{\small }{\bf Abstract.} We study a discrete time multitype branching random walk on a finite space with finite set of types. Particles follow a Markov chain on the spatial space whereas offspring distributions are given by a random field that is fixed throughout the evolution of the particles. Our main interest lies in the averaged (annealed) expectation of the population size, and its long-time asymptotics. We first derive, for fixed time, a formula for the expected population size with fixed offspring distributions, which is reminiscent of a Feynman-Kac formula. We choose Weibull-type distributions with parameter $1/\rho_{ij}$ for the upper tail of the mean number of $j$ type particles produced by an $i$ type particle. We derive the first two terms of the long-time asymptotics, which are written as two coupled variational formulas, and interpret them in terms of the typical behavior of the system. 
\end{quote}

\setcounter{section}{0} 

\setcounter{tocdepth}{2}

%\tableofcontents

%\section{Introduction}\label{sec-Intro}
\section{Introduction.}\label{sec-Intro}

\noindent Branching processes and their applications have been studied for a long time in the communities of mathematics, physics and biology. Additional structure and features like (1) spatiality (random migration of the particles), (2) several types of particles, or (3) dependence of the branching rates on the space and on an additional independent random input lead to variants called (1) {\it branching random walks}, (2) {\it multitype branching processes} and (3) {\it branching processes in random environment}, respectively.  They have attracted researchers because of the additional mathematical richness that they bring into the model and also because of a greater degree of applicability to real-world phenomena. For example, multitype branching processes, where particles are of several types and can give birth to particles of a type other than their own, are motivated by a rich class of applications in physics and biology, such as cosmic ray cascades, bacterial populations and cancer research (see the 
books \cite{M71}, \cite{J75} and \cite{KA02} for various examples).

In this paper, we study a model that exhibits all the three features mentioned above; a {\it multitype branching random walk in random environment}, where particles move independently following a Markov chain on a finite state space and multitype branching takes place according to random site-dependent distributions, which are static, i.e., fixed throughout the evolution. The particles can have a finite number of types, and particles of each type can give birth to new particles of any type. To the best of our knowledge, this combination of features has not yet been considered in the mathematical literature.

One of the most fundamental objects to study is the total number of particles in the system, the global population size, at a given time $n$ in the limit as $n\to\infty$. The two main goals of the present paper are (1) a formula for its expectation (taken over branching/killing and migration, but not over random offspring distributions) for fixed $n$, and (2) the description of the large-$n$ asymptotics of its expectation over the branching probabilities in terms of two coupled variational formulas and their interpretation in terms of a pathwise behavior of the branching process. Hence, we are concerned with the {\it annealed} setting. Note that the corresponding quenched setting is not interesting, since we decided (to keep the technical difficulties low) to work on an arbitrary finite  state space, and therefore there is only a fixed number of random branching probabilities; no effect from infinitely many random inputs is present. 

Fixing $n$, for multitype processes without spatiality, it is standard knowledge \cite{AN72} that the expectation of the population size over the branching/killing for fixed branching probabilities can be described in terms of the $n$-th power of a certain characteristic matrix, the {\it mean-offspring matrix}. For branching random walks, it can be described by the discrete-time version of a Feynman-Kac-type formula. But the question is what structure arises in the case of a multitype branching random walk, i.e., when these two features are combined. Our first result gives the answer: it can be represented in both ways, but on an enlarged space, the product of the state space and the type space. However, there are some subtle differences, one of which is the way in which the branching probabilities are attached to the two spaces: actually, they are attached to the {\it sites} of the state space and to the  {\it edges} of the type space, at least in the particular model that we study.

Turning now to the description of the large-$n$ asymptotics of the annealed expectation of the global population size, we therefore have the choice between working on the expectation of the $n$-th power of a certain random matrix or on the expectation of a Feynman-Kac-like formula. Each of the two settings gives rise to an interesting proof. We decided to work on the latter and to comment on the former only in an informal manner, since it appears to us as if making this line of arguments rigorous would be technically much more involved. We choose the distribution of the branching probabilities as a {\it Weibull-type distribution}, as this distribution exhibits an interesting competition between the migration and the branching probabilities, such that an appealing picture arises. Working on the expectation of the Feynman-Kac-like formula requires the application of a standard large-deviation principle for the empirical pair measures of the underlying Markov chain and goes through quite smoothly. The arising 
two coupled variational formulas give rise to a deeper understanding of the main branching process trees, i.e., of those that give the main contribution to the expected population size. 

Let us give some comments on the existing literature on branching processes in a random environment. One of the models that have been studied before is multitype branching processes with environments varying in time. To name a few, in \cite{T81}, various classification results depending on the long time behavior of the multiplications of mean offspring matrices have been proved, and a much finer analysis in a very general set up has been done in \cite{BCN99} using harmonic functions of multiplication of mean matrices. Let us proceed with spatial branching processes in random environment. Branching discrete-time random walks on $\Z^d$ with time-space i.i.d.~offspring distributions were studied in the context of survival properties, global/local growth rates and diffusivity; and their connections with directed polymers in random environment, see e.g.~\cite{BGK05,Y08,CY11}. Detailed analyses of recurrence/transience properties of discrete-time branching Markov chains with only space-dependent environment, which 
does not exhibit in general the the usual dichotomy valid for irreducible Markov chains, were carried out in \cite{CMP98, MP00, MP03, CP07, M08, BGK09, GMPV10}, to mention some. The main techniques in these studies relate these models to the better-known random walk in random environments, using the spectral properties of underlying Markov process and studying the embedded Galton-Watson processes in random environment.

The remainder of Section~\ref{sec-Intro} is organized as follows. In Section~\ref{Discrete time case} we introduce the branching process in a fixed medium and give representations of the main object of our interest, the global particle number, in Section~\ref{sec-representations}. The random environment is introduced in Section~\ref{sec-randomenvironment}, and our main result is presented in Section~\ref{sec-result}. We comment and interpret it in Section~\ref{sec-explanation}, where we in particular analyse the main quantities appearing in the main result. The special case where migration is dropped leads to even more explicit formulas, which we present in Section~\ref{sec-nomigration}. Finally, in Section~\ref{sec-Frobenius} we phenomenologically discuss and compare another approach to the main result in terms of Frobenius eigenvalue theory for a random matrix, which gives some interesting insights. The proofs of all our results are in Section~\ref{sec-AnneMom}.

\subsection{The Model.}\label{Discrete time case}
 
\noindent The model we will study is a multitype branching Markov chain on a finite state space in discrete time in a fixed environment of branching probabilities. (This environment will be taken random in Section~\ref{sec-randomenvironment} below.) 

Let $(X_n)_{n\in\N_0}$ be an irreducible Markov chain on a finite state space $\mathcal X$ with transition matrix $P=(P_{xy})_{x,y\in\mathcal X}$. Let $\Tcal$ be a finite set, the set of types. We equip $\Tcal$ with a set $\Acal$ of directed edges $(i,j)\in\Tcal\times\Tcal$ and obtain a directed finite graph $\mathcal{G}=(\T,\mathcal{A})$. We assume that each directed edge appears at most once in $\Acal$, and for each $i\in\Tcal$, there is at least one $j\in\Tcal$ such that $(i,j)\in\Acal$. Self-edges $(i,i)$ may appear in $\Acal$. Finally, we assume $\mathcal{G}$ is a connected graph.

To each $y\in\mathcal{X}$ we attach a matrix $F_y=(F_y^{\ssup{i,j}})_{(i,j)\in\mathcal{A}}$ of probability distributions on $\N_0$, the environment. Given $F=(F_y)_{y\in\X}$, we define a discrete-time Markov process $(\eta_n)_{n\in\N_0}$ on $\N_0^{\T\times \X}$, where $\eta_n(i,x)$ is the number of particles of type $i$ at site $x$ at time $n$. The environment $F$ does not depend on time and is fixed throughout the evolution of particles. We specify the transition mechanism of $(\eta_n)_{n\in\N_0}$ as follows: given that the configuration is equal to $\eta$ at time $n$, during the time interval $(n,n+1)$,
\begin{enumerate}
\item a particle of type $i$ located at site $y\in\X$ produces, independently for $j\in\mathcal{T}$ such that $(i,j)\in\A$, precisely $k$ particles of type $j$ at the same site $y$ with probability $F_y^{\ssup{i,j}}(k)$, for any $k\in\N_0$. All offspring productions are independent over all the particles in $\mathcal{X}$ and over time $n\in\N_0$,

\item immediately after creation, each new particle at $x$ chooses a site $y$ with probability $P_{xy}$ and moves there. All jumps are independent over all the particles in $\mathcal{X}$ and over time $n\in\N_0$.
\end{enumerate}

The resulting particle configuration is $\eta_{n+1}$. Note that, unlike in the most general set up of multitype branching processes, we assume that a particle of type $i$ produces particles of type $j$ independently in $j$, that is, the offspring distribution coming from an $i$ type particle is in a product form. Finally, one can immerse the spatial movement into a multitype setting by simply adopting the spatial points as types and rewrite offspring distributions, this time involving terms from the transition matrix $P$. We will comment more on this connection later. 

For definiteness, we consider localized  initial conditions in $\X$ and $\T$. To this end, fix a site $y\in\mathcal{X}$ and type $j\in\T$. We start the Markov chain $(\eta_n)_n$ with the initial configuration $\eta_0(i,x)=\delta_{j}(i)\delta_{y}(x)$, and by ${\tt P}_{j,y}$ and ${\tt E}_{j,y}$ we denote its distribution and expectation, respectively.  Note that they depend on the realization of the environment $F$.

We are interested in the expectation of the global population size, $|\eta_n|:=\sum_{i\in\T,x\in\X}\eta_n(i,x)$, 
\begin{equation}
 u_n(i,x):={\tt E}_{i,x}[|\eta_n|],\qquad n\in\N_0, x\in\mathcal{X}, i\in\Tcal.
\end{equation}
Note that the expectation is taken only on the migration, and the branching/killing, but the environment $F$ is kept fixed.

\subsection{Representations of the expected particle number.}\label{sec-representations}

\noindent Our analysis of the population size is based on a description that is reminiscent of the Feynman-Kac representation of the solution to the heat equation with additive potential. To formulate this, we need to introduce a Markov chain $T=(T_n)_{n\in\N_0}$ on the type space $\Tcal$ with transition probabilities
\begin{equation}\label{pdegdef}
p_{ij}=\frac{\mathds{1}\{(i,j)\in\A\}}{\text{deg}^+(i)},\qquad i,j\in\Tcal,
\end{equation}
where $\text{deg}^+(i)=|\{k\in\Tcal\colon (i,k)\in\A\}|$ is the outdegree of $i$. We define $T$ and $X$ independently on a common probability space and write $\P^{\ssup {T,X}}_{i,x}$ and $\E^{\ssup{T,X}}_{i,x}$ for probability and expectation, respectively, where $T$ starts from $i$ and $X$ from $x$. We denote by $m_{ij}(y)=\sum_{k\in\N_0}kF_y^{\ssup{i,j}}(k)$ the expectation of $F_y^{\ssup{i,j}}$ (the mean number of $j$ type particles at site $y$ that an $i$ type particle at site $y$ produces in one generation) where we set $m_{ij}(y)=0$ for $(i,j)\notin \A$.

The first of the two following representations for $u_n$ is the announced Feynman-Kac-type formulation, which we will use for our proofs in Section~\ref{sec-ProofMain}, and the second one is in terms of the $n$-th power of a particular matrix, which we will use in our heuristic explanations in Section~\ref{sec-Frobenius}. 

\begin{prop}\label{Prop-FKrepr}
For any $i\in\Tcal$ and any $x\in\X$ and any $n\in\N_0$, 
\begin{eqnarray}
u_n(i,x) &=&   \E^{\ssup{T,X}}_{i,x}\Big[  \prod_{l=1}^{n}\big( m_{T_{l-1}T_l}(X_{l-1})\deg^+(T_{l-1}) \big) \Big],\label{AnneMom2} \\
&=&\sum_{j\in\T,y\in\X}B^n_{(i,x),(j,y)},\label{AnneMom2a} 
\end{eqnarray}
where $B^n$ is the $n$th power the $(\Tcal\times\X)\times (\Tcal\times\X)$ matrix $B$ with coefficients
\begin{equation}\label{Bdef}
B_{(i,x),(j,y)}=m_{ij}(x)P_{xy}\mathds{1}\{(i,j)\in\A\}.
\end{equation}
\end{prop}

As mentioned in the previous section, one can enlarge the type space to $\T\times\X$, that is a particle of type $i\in\T$ at site $x\in\X$ can be viewed as a particle of type $(i,x)$. Hence, our model can be seen as a multitype branching process on the finite type space $\T\times \X$. For a general multitype branching process with finitely many types it is very well-known that the expected number of particles at generation $n$ can be described in terms of the entries of the $n$-th power of the mean matrix of one generation. In our model, first a particle produces offsprings and then they migrate along the spatial space. As a result, in the enlarged type space of $\T\times\X$ the mean number of $(j,y)$ offsprings produced by an $(i,x)$ type particle is simply $m_{ij}(x)P_{xy}\mathds{1}\{(i,j)\in\A\}$ which leads to formula in (\ref{AnneMom2a}) with the mean matrix $B$ given as in (\ref{Bdef}). 

An interpretation of formula (\ref{AnneMom2}) is that it is a change of the order of integration. More precisely, one can write the expectation $u_n(i,x)$ as a summation over paths $(x=x_0,i=t_0),(x_1,t_1),\dots,(x_n,t_n)$ in $T\times\X$ of length $n$, the expectation of the number of such paths one can find in the $n$-level branching tree of the process. Clearly, this expectation is given by the term $\prod_{l=1}^n m_{t_{l-1},t_l}(x_{l-1})P_{x_{l-1}x_l}$. Finally, by adding the terms $\deg^+(t_{l-1})$, one can make turn the summation into an expectation over $\P^{\ssup {T,X}}_{i,x}$. 

Writing $\prod_{l=1}^n \cdots$ as $\exp(\sum_{l=1}^n\log\cdots)$ in the expectation in \eqref{AnneMom2}, we encounter a discrete-time version of a Feynman-Kac formula for the Markov chain $(T,X)$ on $\Tcal\times \X$, however with an interesting difference: the potential $\log m_{i,j}(x)$ depends on the {\it vertices} of the space $\X$, but on the {\it edges} of the type space $\Tcal$. A further (more or less negligible) difference is the appearance of the degree term, which accounts for the missing probability structure that we had {\it a priori} on $\Tcal$ and artificially inserted; this term will drop out in the end. Finally, note that Proposition \ref{Prop-FKrepr} reveals that our model can be seen as a version of a branching random walk on the enlarged spatial space $\T\times \X$.

\subsection{The random environment.}\label{sec-randomenvironment}

\noindent Let us describe our assumptions on the random environment. We assume that the collection of all distributions $F_y^{\ssup{i,j}}$ with $y\in\X$ and $(i,j)\in\A$ is independent. Their distribution depends on $(i,j)$, but not on $y$. We call $F=(F_y)_{y\in\X}$ the {\it random environment} and denote by $\Prob$ and $\langle\cdot\rangle$ probability and expectation with respect to $F$, respectively. Note that, $m_{ij}(y)$ is a random variables whose distribution is induced by $F_{y}^{(i,j)}$.  Since we are here interested only in the expectation of the global particle number, we will make our assumptions on the environment only in terms of the quantities $m_{ij}(y)$. In particular, we assume that the collection of the $m_{ij}(y)$ is independent in $y\in\mathcal{X}$ and $i,j\in\Tcal$. Then, it is our goal to find the first two terms of the large-$n$ asymptotics of the expectation of $u_n(i,x)$, i.e., the annealed moments of the global population size when the system starts from one $i$ type particle at 
$x\in\X$.

One can already guess from \eqref{AnneMom2} that the large-time asymptotics of the branching process does not depend on characteristics like expectation or variance of the offspring expectation, but predominantly on their upper tails, since $u_n$ is basically a product of high powers of them. We will study the case where $ m_{ij}(y)$  lies, in terms of upper tails,  in the vicinity of the {\it Weibull distribution} with parameter $1/\rho_{ij}\in(0,\infty)$, i.e.,
    \begin{equation}
    \Prob(m_{ij}(y)>r)\approx\exp\{-r^{1/\rho_{ij}}\}, \qquad r\to\infty. \label{DExpTail}
    \end{equation}
In the language of \cite{GM98}, $\log m_{ij}(y)$ lies in the vicinity of the double-exponential distribution, which is nothing but a reflected Gumbel distribution. The precise assumption on $m_{ij}(y)$ can be written down in terms of the logarithmic moment generating functions given by
\begin{equation}
H_{ij}(t):= \log\langle m_{ij}(y)^t\rangle, \qquad t > 0, \quad i,j \in \mathcal{T}. \label{Hdef}
\end{equation}
Via Tauberian theorems, the upper tails of $m_{ij}(x)$ stand in a one-to-one relation with the regularity of the moment generating function at infinity. Therefore, we sharpen the assumption (\ref{DExpTail}) by requiring the following:

\begin{Assumption}\label{DoubleExp} For any $(i,j) \in \mathcal{A}$ there exists $\rho_{ij}\in(0,\infty)$ such that
    \begin{equation}\label{AssumptionHl}
    \lim_{t\rightarrow\infty}\frac{H_{ij}(ct)-cH_{ij}(t)}{t}=\rho_{ij} c\log c,\qquad c\in(0,1). 
    \end{equation}
\end{Assumption}

For $(i,j)\notin\Acal$, we put $\rho_{ij}=0$. Hence, our environment distribution is characterized by the matrix-valued parameter $\rho=(\rho_{ij})_{i,j\in\Tcal}$. The larger $\rho_{ij}$ is, the thicker the tails of $ m_{ij}(y)$ are, i.e., the easier it is for $m_{ij}(y)$ to achieve extremely high values.

Assumption~\ref{DoubleExp} was used in a number of papers on the parabolic Anderson model on $\Z^d$ (see \cite{GM98}, e.g.) in a similar context. By the virtue of \eqref{AssumptionHl}, in our main result we will see {\it two} explicit terms of the asymptotics of the expectation of $u_n$, both of which describe interesting aspects of the long-time behavior of the branching process. The fact that the parameter $t$ appears both in the arguments of $H$ in the numerator and in the denominator is the reason that this distribution class is particularly amenable to an interesting asymptotics of the annealed asymptotics, since it leads to a match of the large deviation scales of the probability with respect to the Markov chain and the second order term coming from the random medium. For any other regularity assumption on $H$, it would not be possible to match these two scales. This is why we found the Weibull distribution particularly suitable. Of course, this would change if we would work on $\Z^d$ or $\R^d$ instead 
of a fixed finite state space, and many interesting additional questions would arise from the unboundedness of the space, but this is not the focus in the present paper.

\subsection{The main result.}\label{sec-result}

\noindent We are now heading towards a formulation of our main result on the asymptotics of the annealed moments of the global particle number. For any discrete set $S$, we denote by $\Mcal_1(S)$ the set of probability measures on $S$ and by $\Mcal_1^{\ssup {\rm s}}(S^2)$ the set of probability measures on $S^2$ with equal marginals. The first quantity of interest is 
\begin{equation}\label{lambdadef}
 \lambda(\rho) = \sup \Big\{ \langle\mu,\rho\rangle\colon 
\mu \in \mathcal{M}_1^{\ssup {\rm s}}(\mathcal{T}^2) \Big\},\qquad\mbox{where }\langle\mu,\rho\rangle=\sum_{(i,j)\in\A } \mu(i,j)\rho_{ij},
\end{equation}
and the set of the corresponding maximizers:
\begin{equation}\label{Lambdadef}
 \Lambda(\rho):=\Big\{\mu\in\Mcal_1^{\ssup {\rm s}}(\Tcal^2)\colon \langle\mu,\rho\rangle=\lambda(\rho)\Big\}.
\end{equation}
We introduce some notation. Each measure $\nu\in \Mcal_1^{\ssup {\rm s}}((\T\times\X)^2)$ has a number of marginal measures that are defined on different spaces, but in order to keep the notation simple, we denote by $\overline{\nu}$ all these marginals, namely,
\begin{equation}\label{projections}
\begin{aligned}
\overline{\nu}(i,j,x)&=\sum_{y\in\X}\nu((i,x),(j,y)),\qquad\overline{\nu}(i,x)=\sum_{j\in\T}\overline{\nu}(i,j,x),\\
\overline{\nu}(i,j)&=\sum_{x\in\X}\overline{\nu}(i,j,x),\qquad \overline{\nu}(i)=\sum_{j\in\T}\overline{\nu}(i,j).
\end{aligned}
\end{equation}

To describe the second term in the asymptotics, we need to introduce two functionals on measures $\nu \in \Mcal_1^{\ssup s} (( \mathcal{T} \times \X )^2)$, an  energy functional $\Scal$ and an entropy functional $\Ical$. Indeed, define
\begin{eqnarray}
 \Scal(\nu)&:=&\sum_{(i,j)\in\A} \rho_{ij}\sum_{x\in\X}\overline{\nu}(i,j,x)\log \overline{\nu}(i,j,x)+\sum_{(i,j)\in\A}\overline{\nu}(i,j) \rho_{ij}\log \rho_{ij},\label{Sdef}\\
\mathcal{I}(\nu) &:=& \sum_{i,j \in \mathcal{T}}\sum_{x,y \in {\X}} \nu((i,x),(j,y)) \log \frac{\nu((i,x),(j,y))}{\overline{\nu}(i,x) P_{xy} \mathds{1}\{(i,j)\in\Acal\}}.\label{Idef}
    \end{eqnarray}
We set $\mathcal{I}(\nu)=\infty$ if $\nu $ is not absolutely continuous with respect to the measure $((i,x),(j,y))\mapsto \overline{\nu}(i,x) P_{xy}\mathds{1}\{(i,j)\in\A\}$. Then $\mathcal{I}(\nu)$ is equal to the entropy of $\nu$ with respect to this measure; note that it is not normalized, but has mass equal to $\sum_{i\in\Tcal}\bar{\nu}(i)\deg^+(i)$. 

Now we can state our main result: 

\begin{theorem}\label{thm-Main} Under Assumption~\ref{DoubleExp} for any $i\in\T$ and $x\in\X$, as $n\to\infty$,
    \begin{equation}\label{AnneMomAsy}
    \langle u_n(i,x) \rangle = (n!)^{\lambda(\rho)}\e^{-n\chi(\rho)}\e^{o(n)}=  \exp
\Big(\lambda(\rho)n \log \frac{n}{\e}    - n \chi(\rho) + o(n)\Big),
    \end{equation}
where
\begin{equation}\label{Chi}
\chi(\rho) = \inf \Big\{ \mathcal{I}(\nu)-
\Scal(\nu) \colon \nu \in \mathcal{M}_1^{\ssup s}\big( ( \mathcal{T} \times \X )^2 \big), \overline{\nu}\in \Lambda(\rho)\Big\}.
    \end{equation}
\end{theorem}

The proof of Theorem~\ref{thm-Main} is in Section~\ref{sec-ProofMain}. We proceed with some comments on this proof and the interpretation of the formula. Starting from the representation in \eqref{AnneMom2} in Proposition~\ref{Prop-FKrepr}, we follow the patterns of \cite{GM98}, however with some notable changes. {The main step is rewriting the Feynman-Kac representation in terms of the {\it empirical pair measure}
\begin{equation}\label{LocTimes}
     \nu_n = \frac{1}{n} \sum_{l=1}^n \delta_{((T_{l-1},X_{l-1}),(T_l,X_l))},
\end{equation}
which is the central object in this approach. In terms of the space-type random walk $(X,T)$, the number $n \nu_n((i,x),(j,y))$ plays the role of the number of $j$ type offspring of any $i$ type particle located at $x$ by time $n$ that makes a step to $y$ right after creation. Hence, $\nu_n$ stands for the union of all $n$-step paths $((X_0,T_0),\dots,(X_n,T_n))$ that make precisely $n \nu_n((i,x),(j,y))$ steps $(i,x)\to (j,y)$ for every $i,j\in\Tcal$ and every $x,y\in\Xcal$. The term $\Ical(\nu)$ is the negative exponential rate of the probability of this union under the Markov chain $X$, together with the combinatorial complexity of the trajectories of types, and $\Scal(\nu)$, together with the leading term $\lambda(\rho)$, is the one under the expectation w.r.t.~the random environment under Assumption~\ref{DoubleExp}.}

\subsection{Discussion of the variational formulas.}\label{sec-explanation}

\noindent Theorem~\ref{thm-Main} in particular shows that the main contribution to the annealed moments of the particle numbers, $\lambda(\rho)$, comes from those $n$-step branching process subtrees which produce, for some $\mu\in\Lambda(\rho)$, at approximately $n\mu(i,j)$ of the $n$ steps  a number of $j$ type particles from one or more $i$ type particles, for any $i,j\in \Tcal$. Then the value $\langle\mu,\rho\rangle$ gives the leading contribution on the scale $n\log\frac n\e$. It is interesting to note that the optimality of the leading term has nothing to do with the spatial part of the branching process, but only with the creation of particles. The reason is that all the probabilities of spatial actions, i.e., of the random walk $X$, are on the scale $n$, but the values of the offspring expectations $m_{ij}(x)$, averaged over the environment are typically on the scale $n^{O(1)}$ under Assumption~\ref{DoubleExp}.

The interpretion of the second-order term is that, for any maximizer $\nu$ of $\Scal-\Ical$ satisfying $\overline\nu=\mu$, the main contribution comes from those $n$-step branching process trees that place all the births of $j$ type particles from $i$ type particles in such a way on $\Xcal$ that approximately $n\nu((i,x),(j,y))$ such births take place at $x$, and the newly created particle immediately jumps to $y$, for any $i,j\in\Tcal$ and any $x,y\in\X$.

In this light, let us analyse the leading term $\lambda(\rho)$ a bit more closely. A {\it simple cycle} on $\mathcal{G}$ is a path $\gamma=(i_1,\dots,i_l,i_{l+1})$ in $\mathcal{T}$, with steps $(i_{m},i_{m+1})$ in $\Acal$, that begins and ends at the same vertex $i_1=i_{l+1}$, but otherwise has no repeated vertices or edges. We write $(i,j)\in\gamma$ if the directed edge $(i,j)$ belongs to $\gamma$, that is, if $(i,j)=(i_m,i_{m+1})$ for some $m\in\{1,\dots,l\}$. We call $|\gamma|=l$ its {\it length}. We denote by $\Gamma_l$ the set of all simple cycles of length $l$ and by $\Gamma$ the set of all simple cycles. We define
\begin{equation}\label{nugammadef}
\mu_{\gamma}(i,j)=
\begin{cases}1/|\gamma|&\text{if }(i,j)\in \gamma,\\
0&\text{otherwise}.
\end{cases}
\end{equation}
It is clear that $\mu_\gamma\in\Mcal_1^{\ssup{\rm s}}(\Tcal^2)$ for any $\gamma\in \Gamma$. Simple cycles are important for the asymptotics of the annealed moments because the set of extremes of $\Mcal_1^{\ssup {\rm s}}(\T^2)$ consists exactly of the simple cycles of the graph $\mathcal{G}$. Since we could not find a proper reference for this fact, we formulate it as a lemma and prove it in Section ~\ref{sec-ProofMaxCycle}.  

\begin{lemma}\label{extrarecycles}
The set of extremes of the convex set $\Mcal_1^{\ssup {\rm s}}(\T^2)$ is equal to $\{\mu_\gamma\colon\gamma\in\Gamma\}$.
\end{lemma}

Since the optimization problem in (\ref{lambdadef}) is a linear optimization problem on the convex, compact set $\Mcal_1^{\ssup{\rm s}}(\T^2)$, the Krein-Milman theorem and Lemma \ref{extrarecycles} imply the following characterization of the leading term in \eqref{AnneMomAsy}:

\begin{lemma}\label{lem-MaxCyclesLemma}
\begin{equation}\label{MaxCycles}
\lambda(\rho) =  \max\Big\{\langle \mu_\gamma,\rho\rangle;\gamma\in\Gamma\Big\} =\max\Big\{ \frac{1}{|\gamma|} \sum_{m=1}^{|\gamma|} \rho_{i_{m}i_{m+1}}; (i_1,\dots,i_{|\gamma|+1})\in \Gamma \Big\}.
\end{equation}
\end{lemma}

The interpretation of Lemma~\ref{lem-MaxCyclesLemma} is that the leading contribution to the expected population size comes from optimal cycles $(i_1,\dots,i_{|\gamma|})\in\Gamma_{|\gamma|}$ in the sense that already all those $n$-step branching process trees contribute alone optimally which produce only $i_{m+1}$ type particles from $i_m$ type particles for any $m\in\{1,\dots,|\gamma|\}$ (with $i_{|\gamma|+1}=i_1$), but essentially no other offspring.

\subsection{Dropping the migration.}\label{sec-nomigration}

\noindent In this section we illustrate our result in Theorem~\ref{thm-Main} in the special case where migration is absent, i.e., a multitype branching model in a random environment of branching rates without any reference to a spatial component. Here we can give a more explicit description of the  annealed asymptotics. Formally, we drop the spatial component of the model by picking a trivial  Markov chain on $\{x\}$ with $P_{xx}=1$ for some $x$ and remove it from the notation. We analyse the decisive quantities $\lambda(\rho)$, $\Lambda(\rho)$ and $\chi(\rho)$ defined in (\ref{lambdadef}) and (\ref{Lambdadef}). The functionals  $\mathcal S$ and $\mathcal I$ are accordingly modified as follows. For $\nu\in\Mcal_1^{\ssup{\rm s}}(\Tcal^2)$ let
\begin{eqnarray}
\mathcal S(\nu)&:=&\sum_{(i,j)\in\A}\rho_{ij}\nu(i,j)\log \nu(i,j)+\sum_{(i,j)\in\A}\nu(i,j)\rho_{ij}\log\rho_{ij},\\
\mathcal I(\nu)&:=&\sum_{(i,j)\in\A}\nu(i,j)\log\frac{\nu(i,j)}{\overline{\nu}(i)}
\end{eqnarray}
where $\overline{\nu}(i)=\sum_{j}\nu(i,j)$. Let 
\begin{equation}
\Gamma(\rho)=\Big\{\gamma\in \Gamma\colon \mu_{\gamma}\in\Lambda(\rho)\Big\}=\Big\{\gamma\in\Gamma\colon \langle \mu_\gamma,\rho\rangle=\lambda(\rho)\Big\}.
\end{equation}

We restrict to the case where $\rho_{ij}=\rho_i\geq 1$ for all $i,j$. We give an explicit solution to the variational formula in \eqref{Chi} that defines $\chi(\rho)$. In the even more restricted case where $\rho_{ij}=\rho\in(0,\infty)$ for all $i,j$, we have $\lambda(\rho)=\rho$ and $\Lambda(\rho)=\Mcal_1^{\ssup{\rm s}}(\T^2)$. Let $l_{\min}$ be the girth of the directed graph $\mathcal{G}$, that is, the length of a shortest simple cycle of the graph.

\begin{lemma}\label{Mut-prop2}
Let $\rho_{ij}=\rho_i\geq 1$ for all $j$. Then 
\begin{equation}\label{newvar}
\chi(\rho)=\min\Big\{\lambda(\rho)\log |\gamma|-\sum_{m=1}^{|\gamma|}\frac{1}{|\gamma|}\rho_{i_m}\log \rho_{i_m}\colon \gamma\in\Gamma(\rho) \Big\}.
\end{equation}
More specifically, if $\rho_i= \rho\geq 1$ for all $i$, then $\chi(\rho)=\rho\log l_{\min}-\rho\log\rho$, and the set of minimizers in (\ref{newvar}) is equal to $\Gamma_{l_{\min}}$.
\end{lemma}

\subsection{Comparison to a Frobenius eigenvalue approach.}\label{sec-Frobenius}

\noindent To analyse the annealed moment asymptotics of the particle number in the spatial multitype branching process in random environment, another approach is also very tempting.  We want to roughly explain this briefly for the simpler case where migration is absent, i.e., in the setting of Section~\ref{sec-nomigration}. (The general case is not much different, thanks to \eqref{AnneMom2a}, but notationally more cumbersome.) In this case, \eqref{AnneMom2} (or Lemma~\ref{thm-MomFormLemma} below) simplifies to 
\begin{equation}\label{AnneMom3}
u_n(i) = \sum_{j\in\T}\big(M^n\big)_{ij},\qquad i\in\Tcal, n\in\N_0,  
\end{equation}
where we recall that the left-hand side is the expected number (taken only over the branching/killing) of particles at time $n$ when we started with just one $i$ type at time $0$, for fixed environment, and we write $M=(m_{ij})_{i,j\in\Tcal}$ for the expectation matrix, where $m_{ij}$ is the expectation of the number of $j$ type offspring of a $i$ type particle. For simplicity, we assume that $\Acal=\Tcal^2$, i.e., that Assumption~\ref{DoubleExp} holds for every $(i,j)\in\Tcal^2$. Assume that, almost surely, $M$ is irreducible, that is, for any $i,j\in\Tcal$ there is an $n\in\N$ such that $(M^n)_{ij}>0$. Then the high powers of $M$ can be approached with the help of the {\it Frobenius eigenvalue} as follows. Define, for any positive irreducible matrix $A=(a_{i,j})_{i,j\in\Tcal}$,
\begin{equation}\label{Frobeniusdef}
\mu(A)=\lim_{k\to\infty}\frac 1k\log \sum_{j}(A^k)_{i,j},
\end{equation}
and note that this limit exists and does not depend on $i$ \cite{S06}. Furthermore, $\e^{\mu(A)}$ is equal to the largest eigenvalue of $A$, the Frobenius eigenvalue of $A$, which has also the characteristic property that it is a simple eigenvalue both algebraically and geometrically. Hence, we can approximate $\langle u_n(i)\rangle\approx \langle \e^{n\mu(M)}\rangle$. This means that we are faced with the question of a large-deviation principle for the random vector $M$. 

We consult Assumption~\ref{DoubleExp} and see that, for each $(i,j)\in\Tcal^2$, the variable $\log m_{ij}-H_{ij}(n)/n$ satisfies an LDP with speed $n$ and rate function $\R\ni m\mapsto \frac{\rho_{ij}}\e\e^{m/\rho_{ij}}$. This is easily calculated from \eqref{DExpTail} by using that $H_{ij}(n)=\rho_{ij} n\log(\rho_{ij} n)-\rho_{ij} n+o(n)$ for $n\to\infty$. By independence, the entire matrix $K_n=\log M-H(n)/n=(\log m_{ij}-H_{ij}(n)/n)_{i,j\in\Tcal} $ satisfies an LDP on $\R^{\Tcal^2}$ with speed $n$ and rate function $m\mapsto \prod_{i,j}\frac{\rho_{ij}}\e\e^{m_{ij}/\rho_{ij}}$. We have $m_{ij}=\e^{K_n(i,j)}\e^{H_{ij}(n)/n}$ for each $i,j\in\Tcal$. In order to identify the asymptotics of $\langle \e^{n\mu(M)}\rangle$, one needs to employ the LDP that $K_n$ satisfies, but it appears difficult to write the Frobenius eigenvalue of $M$ in terms of the one of $K_n$. For making decisive progress here, it seems as if one must apply some tools employed in the present paper to the formula \eqref{Frobeniusdef} as follows: Explicitly write out the $n$-fold product of the matrix $A$, write it in terms of the empirical pair measure of the resulting multi-index and use the well-known combinatorics for the number of strings that lead to a given empirical pair measure, to get that
\begin{equation}
\mu(A)=\sup_{\nu\in\Mcal_1^{\ssup{\rm s}}(\Tcal^2)}\Big(\langle \nu,\log A\rangle -\sum_{i,j\in\Tcal}\nu(i,j)\log\frac{\nu(i,j)}{\overline\nu(i)}\Big),
\end{equation}
where $\langle \nu,\log A\rangle=\sum_{i,j\in\Tcal}\nu(i,j)\log a_{ij}$, and $\overline\nu$ denotes the marginal measure of $\nu$. It appears also not easy to extract the precise dependence of the leading term $\mu((\rho n/\e)^\rho)$ on $n$ without this procedure. Summarizing, we believe that the method that we use in the present paper is essentially the only doable way.

\section{Proofs of the Main Results.}\label{sec-AnneMom}

\subsection{Proof of Proposition~\ref{Prop-FKrepr}.}\label{sec-ProofProp}

\noindent We now write the expected number of offspring in terms of an expectation of the product of the expectation matrices $M(x)=(m_{ij}(x))_{i,j\in\Tcal}$ along the path of the Markov chain  $X=(X_n)_{n \in \N_0}$ on $\X$. This may be seen as a discrete-time version of a Feynman-Kac formula and is completely standard. By $\P_x^{\ssup X}$ and $\E_x^{\ssup X}$ we denote probability and expectation with respect to the random walk when started at $x\in\X$.

\begin{lemma}\label{thm-MomFormLemma} For any $i,j\in\T$, $x,y\in\X$ and any time $n$,
\begin{equation}\label{ExpValDiscrete}
 {\tt E}_{(i,x)}[\eta_{n}(j,y)]= \E_{x}^{\ssup X}\big[\big( M(X_0)\cdots M(X_{n-1})\big)_{ij} \delta_y(X_n)].
\end{equation}
\end{lemma}

\begin{proofsect}{Proof}
We proceed the proof by induction on $n$. For $n = 1$ it is easy to see that
\begin{equation}
{\tt E}_{(i,x)}[\eta_{1}(j,y)] =   m_{ij}(x)P_{xy}, 
\end{equation}
from which immediately the assertion for $n=1$ follows. Now assume that (\ref{ExpValDiscrete}) holds for $n \in \N$. We have 
\begin{align}
{\tt E}_{(i,x)}[\eta_{n+1}(j,y)]&=\sum_{k\in\T}\sum_{z\in\X}m_{kj}(z)P_{zy}{\tt E}_{(i,x)}[\eta_n(k,z)]
\\&=\sum_{k\in\T}\sum_{z\in\X}m_{kj}(z)P_{zy}\E_{x}^{\ssup X}\big[\big( M(X_0)\cdots M(X_{n-1})\big)_{ik} \,\delta_{X_n}(z))]
\\&=\E_{x}^{\ssup X}\big[\big( M(X_0)\cdots M(X_{n})\big)_{ij} \,\delta_{X_{n+1}}(y))]
\end{align}
where in the second equality we have used the induction step and in the third equality we have used the Markov property of $X_n$. Hence, by induction the proof of the lemma is finished. 
\end{proofsect}

Now, explicitly writing out the matrix product in (\ref{ExpValDiscrete}) and summing over $j\in\T$, we get
\begin{equation}\label{AnneMom1}
u_n(i,x) = \sum_{k_1,\ldots,k_{n} \in \mathcal{T}} \E_{x}^{\ssup X} \Big[  \prod_{l=1}^{n}m_{k_{l-1}k_l}(X_{l-1})  \Big],\qquad k_0=i.   
\end{equation}
Now we absorb the transition probabilities given in \eqref{pdegdef} in the product and easily rewrite the right-hand side of \eqref{AnneMom1} as the right-hand side of \eqref{AnneMom2}, which finishes the proof of Proposition~\ref{Prop-FKrepr}.

\subsection{Proof of Theorem~\ref{thm-Main}.}\label{sec-ProofMain}

\noindent We start from Proposition~\ref{Prop-FKrepr} and recall the normalized empirical pair measure of $(T,X)$:
$$
     \nu_n = \frac{1}{n} \sum_{l=1}^n \delta_{((T_{l-1},X_{l-1}),(T_l,X_l))},
$$
which is a probability measure on $( \mathcal{T} \times \X )^2$. Now fix $i\in\Tcal$ and $x\in\X$ and rewrite the right-hand side of (\ref{AnneMom2}) in terms of $\nu_n$:
\begin{equation}\label{AnneMom4}
 u_n(i,x) =   \E^{\ssup{T,X}}_{i,x}\Big[ \Big( \prod_{k,j \in \mathcal{T}}\prod_{z \in \X} m_{kj}(z)^{n \sum_{ y \in \X} \nu_n((k,z),( j, y))}\Big)\\
  \Big( \prod_{ k\in\T}\text{deg}^+(k)^{n\sum_{ j\in\T,z, y\in\X}\nu_n((k,z),( j, y)) }\Big) \Big].    
\end{equation}
Let us introduce the function $\mathcal{D}$ for $\nu\in\Mcal_1((\T\times \X)^2)$ by
\begin{equation}
\mathcal{D}(\nu):=\sum_{k\in\T}\overline{\nu}(k)\log (\text{deg}^+(k)).
\end{equation}
Now we take expectation with respect to the environment of the right-hand side of (\ref{AnneMom4}) and use that the $m_{kj}(z)$  with $z\in\X$ and $k, j\in\Tcal$ are independent. Using the notation introduced in (\ref{projections}) and the logarithmic moment generating function $H_{k j}(t) = \log\langle m_{k j}(0)^t\rangle$ we get:
\begin{equation} \label{AnneMom5}
\begin{aligned}
\langle u_n(i,x) \rangle &= \Big\langle\E_{i,x}^{\ssup{T,X}}\Big[\exp\Big(\sum_{k,j,z}\big(\log m_{kj}(z)\big)n\bar{\nu}_n(k,j,z)+\sum_k \log(\text{deg}^+(k)) n\bar{\nu}_n(k)\Big)\Big]\Big\rangle \\& =   \E_{i,x}^{\ssup{T,X}}\Big[ \e^{n\mathcal{D}(\nu_n)} \prod_{k, j \in \mathcal{T}} \e^{\sum_{z \in \X} H_{k j}(n \overline{\nu}_n(k, j,z) )} \  \Big] \\
&  =   \E^{\ssup{T,X}}_{i,x} \Big[ \e^{n\mathcal{D}(\nu_n)}\prod_{k, j \in \mathcal{T}} \exp\Big(n\sum_{z \in\X} \frac{H_{k j}(n \overline{\nu}_n(k, j,z) ) -  \overline{\nu}_n(k,j,z) H_{k j}(n)  }  {n}  \Big) \\
& \qquad \times  \exp \Big( \sum_{k, j \in \mathcal{T}}\overline{\nu}_n(k, j) H_{k j}(n)\Big)\  \Big].  
\end{aligned}
\end{equation}
Recall Assumption (\ref{DoubleExp}) and note that $H_{kj}$ satisfies the asymptotics $H_{kj}(n) = \rho_{kj}n\log(\rho_{kj}n) - \rho_{kj}n + o(n)$ as $n \rightarrow \infty$ \cite{GM98}. Furthermore, we use the asymptotics in \eqref{AssumptionHl} for every $z\in \X$ to conclude that as $n\to\infty$
\begin{equation}\label{AnneMom6}
\begin{aligned}
\langle u_n(i,x) \rangle & =  \E^{\ssup{T,X}}_{i,x}  \Big[\e^{n\mathcal{D}(\nu_n)} \exp\Big( n\sum_{k, j \in \mathcal{T}} \rho_{k j}\sum_{z \in \X} \overline{\nu}_n(k,j,z) \log \overline{\nu}_n(k, j,z)  \\ 
&\qquad +  \sum_{k, j \in \mathcal{T}} \overline{\nu}_n(k, j) \Big( \rho_{k, j}n \log (\rho_{k j}n) - \rho_{k j}n \Big)\Big)\Big]\e^{o(n)} \\
& =   \E^{\ssup{T,X}}_{i,x} \Big[\e^{ n \Scal(\nu_n)+n\mathcal{D}(\nu_n)} \exp\Big( n \log \Big(\frac n \e\Big) \langle\overline{\nu}_n, \rho\rangle \Big)\Big] \e^{o(n)} ,
\end{aligned}
\end{equation}
where we used the definition of the functional $\Scal$ in \eqref{Sdef} and the notation from \eqref{lambdadef} for the last equality. Also note that in the above display we have $\langle\bar{\nu},\rho\rangle=\sum_{i,j}\overline{\nu}(i,j)\rho_{ij}$. 

Our main tool in the proof is that the pair empirical measure $(\nu_n)_n$ of the Markov chain $(T_n,X_n)_{n\in\N}$ satisfies a large deviation principle on the set of probability measures on $(\Tcal\times\X)^2$ with scale $n$ and the good rate function $I'$ given by 
\begin{equation}
I'(\nu)=\sum_{(i,j)\in\A,x,y\in\X}\nu((i,x),(j,y))\log\frac{\nu((i,x),(j,y))}{\overline{\nu}(i,x)P_{xy}/\text{deg}^+(i)}=I(\nu)+\mathcal{D}(\nu)
\end{equation}
if $\nu$ satisfies the marginal condition (i.e., lies in $\Mcal_1^{\ssup{\rm s}}((\Tcal\times\X)^2)$), and $I'(\nu)=\infty$ otherwise.

\noindent{\underline{\bf Upper bound}:} Since we are on a finite space, it is clear that we can restrict ourselves to $\nu\in\mms$. We define
\begin{equation}
\Lambda_\delta(\rho):=\Big\{\mu\in \mathcal{M}_1^{\ssup{\rm s}}(\T^2)\colon \langle \mu, \rho\rangle\geq \lambda(\rho)-\delta\Big\},\qquad \delta>0.
\end{equation}
We split the expectation on the right-hand side of \eqref{AnneMom6} into the contribution from the events $\{\overline\nu_n\in \Lambda_\delta(\rho)\}$ and its complement. For the first order term in the first part we use the upper bound $\lambda(\rho)$ and, in the second part the first order term is bounded above by $\e^{n\log(n/\e)(\lambda(\rho)-\delta/2)}$. Hence, we obtain, for all sufficiently large $n$,
\begin{equation}
\begin{aligned}
\langle u_n(i,x) \rangle \leq\e^{n\log(n/\e)\lambda(\rho)}\E^{\ssup{T,X}}_{i,x} \Big[ \e^{ n \Scal(\nu_n)+n\mathcal{D}(\nu_n)}\1\big\{\overline{\nu}_n\in\Lambda_\delta(\rho)\big\}\Big]\e^{o(n)}+\e^{n\log(n/\e)(\lambda(\rho)-\delta/2)} .
\end{aligned}
\end{equation}
Since $\Lambda_\delta(\rho)$ is closed and $\Scal$ and $\mathcal{D}$ are continuous and bounded, with the same reasoning as in the proof of the upper bound in Varadhan's Lemma (see e.g. pages 32-34 in \cite{H00}) , we can conclude that as $n\to\infty$
\begin{equation}
\E^{\ssup{T,X}}_{i,x} \Big[ \e^{ n \Scal(\nu_n)+n\mathcal{D}(\nu_n)}\1\big\{\overline{\nu}_n\in\Lambda_\delta(\rho)\big\}\Big]\leq\e^{o(n)}\exp\Big(-n\inf_{\Lambda_\delta(\rho)}(I-\Scal)\Big).
\end{equation}
In the above we also used that $I'=I+D$. Since $I$ and $\Scal$ are continuous and bounded functionals, $\Lambda_{\delta}(\rho)$ is compact for any $\delta\geq 0$ and $\Lambda_{\delta}(\rho)\downarrow \Lambda(\rho)$ as $\delta\downarrow \rho$, we can conclude that
\begin{equation}
 \lim_{\delta\downarrow 0} \inf_{\Lambda_\delta(\rho)}(I-\Scal)\geq \inf_{\Lambda(\rho)}(I-\Scal).
\end{equation}
This finishes the proof of the upper bound in \eqref{AnneMomAsy}.

\noindent{\underline{\bf Lower Bound:}} Now we prove the lower bound in \eqref{AnneMomAsy}. For any $\nu\in\Lambda(\rho)$ let
\begin{equation}
B_{n}(\nu)=\{\nu'\in \mms:d(\nu',\nu)<n^{-1/2}\}
\end{equation}
where $d(\cdot,\cdot)$ is the usual total variation distance on probability distributions. We can find a constant $c\in\R$ that depends only on $\rho$ such that for any $\nu\in\Lambda(\rho)$, for any $n$ and $\nu'\in B_{n}(\nu)$ 
\begin{equation}
\langle\nu',\rho\rangle\leq \lambda(\rho)+ cn^{-1/2}.
\end{equation}
Then, for any $\nu\in\Lambda(\rho)$, from (\ref{AnneMom6}) we have the lower bound
\begin{equation}
\langle u_n(i,x) \rangle \geq \e^{o(n)}\exp\Big({n \log \Big(\frac n \e\Big)(\lambda(\rho)-cn^{-1/2})}\Big)\E^{\ssup{T,X}}_{i,x} \Big[ \1\{\nu_n\in B_n(\nu)\}\e^{ n \Scal(\nu_n)+n\mathcal{D}(\nu_n)}\Big].
\end{equation}
{Using the explicit combinatorics (see the proof of Theorem II.8 in \cite{H00}) we have for any $\nu\in\mms$
\begin{equation}
\liminf_{n\to\infty}\frac{1}{n}\log \P^{\ssup{T,X}}_{i,x}\big(\nu_n\in B_n(\nu)\big)\geq -I'(\nu).
\end{equation}
}
Using the continuity of $\Scal$ and $\mathcal{D}$ we get for any $\nu\in \Lambda(\rho)$
\begin{equation}
 \langle u_n(i,x) \rangle \geq \e^{o(n)}\exp\Big({n \log \Big(\frac n \e\Big)\lambda(\rho)}\Big) \e^{n\Scal(\nu)+n\mathcal{D}(\nu)}\e^{-nI'(\nu)}.
\end{equation}
Recalling that $I'=I+D$ and taking the supremum over $\Lambda(\rho)$ finishes the proof of the lower bound.

\subsection{Proof of Lemma~\ref{extrarecycles}.}\label{sec-ProofMaxCycle} 

\noindent We first prove that for any $\gamma\in\Gamma$ the measure $\mu_\gamma$ is an extremal point of $\Mcal_1^{\ssup{\rm s}}(\T^2)$. We need to show that if $\mu_{\gamma}=a\mu_1+(1-a)\mu_2$ for some $a\in(0,1)$ and $\mu_1,\mu_2\in\ms$, then it must be the case that $\mu_1=\mu_2$.  Since $a\in(0,1)$, both supports of $\mu_1$ and $\mu_2$ are contained in the support of $\mu_\gamma$. Without loss of generality, suppose that there exists an edge $(i,j)$ such that $ (i,j)\in \supp{\mu_\gamma}$ with $\mu_1(i,j)=0$. Since $\gamma$ is a simple cycle, $\supp(\mu_1)\subsetneq \supp{\mu_\gamma}$ and $\mu_1\in\ms$, we get that $\mu_1\equiv 0$. Thus, the supports of $\mu_1$ and $\mu_2$ are both equal to the support of $\mu_\gamma$, which is equal to the set of the edges of the cycle $\gamma$. The only measure in $\ms$ having this support is $\mu_\gamma$. Hence, $\mu_\gamma=\mu_1=\mu_2$. 

Now we show that if $\mu \in \mathcal{M}_1^{\ssup {\rm s}}(\mathcal{T}^2)$ is not of form $\mu_\gamma$ for some $\gamma\in\Gamma$, then it is not an extremal point of $\ms$. By the marginal property, $\supp(\mu)$ contains some simple cycle $\gamma_1$. On the other hand, $\supp(\mu)$ is not a simple cycle. As a result, $\supp(\mu)$ contains a second simple cycle $\gamma_2\in\Gamma$ such that $\gamma_2\not=\gamma_1$. Now we choose $\eps>0$ small enough so that $\mu(i,j)\geq 4\eps$ for any $ (i,j)\in\gamma_1\cup\gamma_2$. Then we have $\mu-\eps|\gamma_1|\mu_{\gamma_1}-\eps|\gamma_2|\mu_{\gamma_2}\geq 2\eps$ and $1-\eps|\gamma_1|-\eps |\gamma_2|\in(1/2,1)$. Then the measure 
\begin{equation}
\mu':=\frac{1}{1-\eps |\gamma_1|-\eps|\gamma_2|}\Big(\mu-\eps|\gamma_1|\mu_{\gamma_1}-\eps|\gamma_2|\mu_{\gamma_2}\Big)
\end{equation}
belongs to $\Mcal_1(\T^2)$. Also since $\mu,\mu_{\gamma_1},\mu_{\gamma_2}$ satisfy marginal property, so does $\mu'$ and hence, $\mu'\in\ms$. Finally, we can write $\mu$ as
\begin{equation}
\mu=\eps|\gamma_1|\mu_{\gamma_1}+\eps|\gamma_2|\mu_{\gamma_2}+(1-\eps|\gamma_1|-\eps|\gamma_2|)\mu'
\end{equation}
and, as noted above $\eps|\gamma_1|,\eps|\gamma_2|,1-\eps|\gamma_1|-\eps|\gamma_2|\in(0,1)$ and $\mu_{\gamma_1},\mu_{\gamma_2},\mu'\in\ms$. Therefore, $\mu$ is not an extremal point of $\ms$.

\subsection{Proof of Lemma~\ref{Mut-prop2}.}\label{sec-ProofMut2}

\noindent We rewrite $\Ical-\Scal=\widetilde G+\widetilde F$, where
\begin{equation}
\widetilde{G}(\nu):=\sum_{i,j\in\T}(1-\rho_i)\nu(i,j)\log\frac{\nu(i,j)}{\overline{\nu}(i)},\qquad \widetilde{F}(\nu):=-\sum_{i\in\T}\rho_i\overline{\nu}(i)\log \overline{\nu}(i)-\sum_{i\in\T}\overline{\nu}(i)\rho_i\log \rho_i.
\end{equation}
Using this and Lemmas \ref{extrarecycles} and \ref{lem-MaxCyclesLemma}, we get that $\Lambda(\rho)$ is the convex hull of $\{\mu_\gamma:\gamma\in\Gamma(\rho)\}$. This and the fact that $\widetilde{F}$ is strictly concave gives
\begin{equation}\label{conbound}
\widetilde{F}(\nu)\geq \min\Big\{\widetilde{F}(\mu_{\gamma})\colon\gamma\in\Gamma(\rho)\Big\}.
\end{equation}
We have $\widetilde{G}(\nu)\geq 0$ for any $\nu\in\Mcal_1((\T\times\X)^2)$, since $\rho_i\geq 1$ for all $i\in\T$. Moreover, for any cycle $\gamma$, $\widetilde{G}(\nu_\gamma)=0$. Hence, by (\ref{conbound}) we get that 
\begin{equation}
\chi(\rho)=\min\Big\{\widetilde{F}(\mu_{\gamma})\colon\gamma\in\Gamma(\rho)\Big\}.
\end{equation}
Finally, note that for any simple cycle $\gamma=(i_1,\dots,i_{|\gamma|})$ with $\gamma\in\Gamma(\rho)$
\begin{equation}
\widetilde{F}(\mu_\gamma)=\lambda(\rho)\log |\gamma|-\sum_{m=1}^{\gamma}\frac{1}{|\gamma|}\rho_{i_m}\log \rho_{i_m}.
\end{equation}
This finishes the first part of the proof. If $\rho_i=\rho$ for all $i$, then $\lambda(\rho)=\rho$, $\Gamma(\rho)=\Gamma$ and $\widetilde{F}(\mu_\gamma)=\rho \log|\gamma|-\rho\log\rho$. Finally, noting that $\widetilde{F}$ is strictly concave finishes the proof of the second part.

%%%%%%%%%%%%%%%%%%%%%%%%%%%%%%%%%%%%%%%%%%%%%%%%

\end{document}